\documentclass[3p, preprint]{elsarticle}
\pdfoutput=1
\usepackage{amsmath}
\usepackage{moreverb}
\usepackage{amssymb}
\usepackage{tabularx}
\usepackage{booktabs}
\usepackage[percent]{overpic}
\usepackage{relsize}
\usepackage{mathtools}
\newcommand{\Ptil}{{\tilde{P}}}

\usepackage{hyperref}

\journal{ArXiv}

\begin{document}

\begin{frontmatter}
\title{Modeling metabolic networks including gene expression and uncertainties\tnoteref{cgs}}

\author[ovgu]{H.~Lindhorst\corref{cor1}}
\ead{henning.lindhorst@ovgu.de}
\author[tub]{S.~Lucia}
\ead{sergio.lucia@tuberlin.de}
\author[ovgu]{R.~Findeisen}
\ead{rolf.findeisen@ovgu.de}
\author[ku]{S.~Waldherr}
\ead{steffen.waldherr@kuleuven.be}

\cortext[cor1]{Corresponding author.}
\address[ovgu]{Otto-von-Guericke-University Magdeburg, Institute for Automation Engineering, Universit\"atsplatz 2, D-39112 Magdeburg, Germany}
\address[tub]{Technische Universit\"at Berlin, Laboratory for Internet of Things and Smart Buildings. Ernst-Reuter-Platz 7, D-10587 Berlin, Germany}
\address[ku]{KU Leuven, Department of Chemical Engineering, Celestijnenlaan 200F bus 2424, 3001 Heverlee, Belgium}

\tnotetext[cgs]{Sponsored by German Federal Research Ministry (BMBF) Fkz.\ 031L0017A}

\begin{abstract}
Constraint based methods, such as the Flux Balance Analysis, are widely used to model cellular growth processes without relying on extensive information on the regulatory features.
The regulation is instead substituted by an optimization problem usually aiming at maximal biomass accumulation.
A recent extension to these methods called the dynamic enzyme-cost Flux Balance Analysis (deFBA) is a fully dynamic modeling method allowing for the prediction of necessary enzyme levels under changing environmental conditions.
However, this method was designed for deterministic settings in which all dynamics, parameters, etc. are exactly known.
In this work, we present a theoretical framework extending the deFBA to handle uncertainties and provide a robust solution.
We use the ideas from multi-stage nonlinear Model Predictive Control (MPC) and its feature to represent the evolution of uncertainties by an exponentially growing scenario tree.
While this representation is able to construct a deterministic optimization problem in the presence of uncertainties, the computational cost also increases exponentially.
We counter this by using a receding prediction horizon and reshape the standard deFBA to the short-time deFBA (sdeFBA).
This leads us, along with further simplification of the scenario tree, to the robust deFBA (rdeFBA).
This framework is capable of handling the uncertainties in the model itself as well as uncertainties experienced by the modeled system.
We applied these algorithms to two case-studies: a minimal enzymatic nutrient uptake network, and the abstraction of the core metabolic process in bacteria.
\end{abstract}


\end{frontmatter}

\section{Introduction}
Bacteria are encountering variations in their natural living spaces and have developed complex regulatory mechanism to cope with these.
Without full knowledge on these mechanisms we usually rely on optimization based methods to predict the behavior of the bacteria.
The most prominent approach utilizing this concept is the Flux Balance Analysis (FBA) \cite{varma1994metabolic}.
This idea optimizes the fluxes in steady state, such that a biomass producing flux is maximized.
This simple approach has led to a number of different extensions like the iterative FBA methods presented in \cite{varma1994stoichiometric} and \cite{mahadevan2002dynamic}.
They calculate optimal fluxes for each iterative step with the standard FBA and update the environmental conditions with regards to the calculated fluxes.
While these models can already include changes in the environment, they tend to react very abruptly to these changes as adaptation processes are not modeled.
These adaptation processes are mainly changes in the utilized pathways, which need the presence of the corresponding enzymes to work correctly.
The relation between the pathway usage and the enzyme levels was first modeled in the Resource Balance Analysis (RBA) \cite{goelzer2011cell}.
These models optimize the growth rate under fixed environment while including enzymatic flux constraints to limit uptake and metabolic fluxes.
The solution are the flux rates and the enzyme distribution under the quasi-stationary assumption, leading to the optimized exponential growth.
These methods culminated in the dynamic enzyme-cost FBA (deFBA) \cite{waldherr2015}, which predicts all reaction rates and the enzyme levels over a time course.
It is especially useful in investigating the adaptation process over the course of nutrient changes, e.g. the switch from aerobic to anaerobic growth conditions or the depletion of a preferred carbon source.

All of these modeling methods rely heavily on experimental measurements for the stoichiometry, flux bounds, catalytic constants, etc.
These measurements naturally include errors and thus introduce a form of uncertainty to the biological models.
While there are some extensions to the FBA handling uncertainties in the stoichiometry of the network \cite{zavlanos2011robust} or the steady-state assumption \cite{almass2014robust}, these are only applicable to static problems without the inclusion of enzyme catalysis.

In this work we present a general framework for the class of FBA methods capable of portraying uncertain parameters, and to ensure robust results with respect to these uncertainties.
It is important to ensure robustness, because otherwise we might get infeasible results, e.g. in a model with uncertain nutrient dynamics a non robust solution might predict uptake of an already depleted source.
With the deFBA being the most complex FBA method to date, we have chosen it as the basis for our framework and combine it with some ideas originating in robust Model Predictive Control (MPC) \cite{bemporad1999robust}.
As most MPC applications naturally include uncertainties, these methods have already shown they can be used in real world applications.
The most classic approach in robust MPC is the min-max MPC \cite{campo1987robust}. 
The solution for this problem is a series of control inputs minimizing the objective for a worst-case realization of the uncertainties.
This approach could be included in the deFBA, but for more complex networks the problem may lead to a very conservative solution or even become infeasible as shown in \cite{scokaert1998min}.
Another possible class of robust methods is the tube-based MPC \cite{rakovic2011fully}.
Here the solution of the nominal control problem is combined with a so called ancillary-controller ensuring the evolution of the real uncertain system is constrained to a predefined tube centered around the nominal trajectory.
This method guarantees stability and the solution satisfies all constraints under the uncertainties.
Yet, optimality of the solution is not addressed, and can therefore not be used in the biological setting, as the optimality of the solution is the governing assumption for all FBA methods.
Instead we use the multi-stage nonlinear MPC (NMPC) \cite{lucia2013}, \cite{lucia2014} as the robust method to be combined with the deFBA.
The multi-stage NMPC uses a representation of uncertainties by means of a scenario tree on a finite prediction horizon.
These trees are constructed by modeling the influence of an uncertainty as a discrete event, which creates additional vertices in the tree.
Hence, the number of paths in this tree grows exponentially in the number of uncertainties and the observed time frame.
This leads to a very large optimization problems if applied to deFBA models with a large number of measured values and a prediction over long time spans.
Consequently, we apply the idea of the receding prediction horizon to the deFBA to minimize the size of the scenario tree and thus the problem size.
Our first contribution is the \emph{short-term deFBA} (sdeFBA), which utilizes the receding prediction horizon to form an iterative scheme to solve deFBA problems.
In this context, we suggest a systematic way to choose the prediction horizon for the sdeFBA such that the results are qualitatively similar or superior to the deFBA solutions.
Using the sdeFBA in itself can be beneficial for some applications, but most importantly is can be used as a predictor inside a MPC scheme.

The next step in formulating the robust framework is the integration of the scenario tree in the sdeFBA.
We call the result the \emph{robust deFBA} (rdeFBA).
While capable to handle all kinds of uncertainties relevant to the metabolic networks, we focus on measurement errors in the catalytic constants.
The effects of uncertainty in these constants are usually underestimated, but should be given proper thought as our examples show.
We consider a minimal example to clarify the differences in the presented method and a medium sized example describing the core metabolic process in bacteria.

In summary, our proposed framework is an extension of the deFBA, allowing us to include different kinds of uncertainties within the model of metabolic networks.
The result follows a strict optimality principle and is guaranteed to be admissible under the influence of the constraints.
Furthermore, we depict the impact of measurement errors by comparison of robust and regular solutions of the examples mentioned above.

\section{Dynamic Enzyme-cost FBA}
Let us recap the deFBA and its mathematical description \cite{waldherr2015}.
This method was developed to analyze the biomass fluxes under changing environmental conditions as well as the enzyme levels necessary to realize these. 

\subsection{Modeling metabolic networks}
The networks we are investigating consist of a set of molecular species, which we divide into external species $Y \in \mathbb{R}^{n_y}$, internal species $X \in \mathbb{R}^{n_x}$, and a set of macromolecules  $P \in \mathbb{R}^{n_p}$.
All species have units of molar amounts, $\left[ X \right] = \left[Y\right] = \left[P \right] = 1~\mathrm{mol}$.
We assign each macromolecule in $P$ its \emph{molecular weight} $b_i,~i \in \{1, \ldots, n_p\}$, which we use to measure the \emph{total biomass} $b^T P,~b=(b_1, \ldots, b_{n_p})$.
We denote the $n_r$ reactions between the species by the resulting reaction fluxes $V = (V_y^T, V_x^T, V_p^T)^T,~\left[V \right]= 1~\frac{\mathrm{mol}}{\mathrm{s}}$, which are divided into three classes:
\begin{itemize}
 \item Exchange reactions $V_y \in \mathbb{R}^{r_y}$ transport matter between the inside and the outside of the cell;
 \item Metabolic reactions $V_x \in \mathbb{R}^{r_x}$ convert the metabolites among each other;
 \item Biomass reactions $V_p \in \mathbb{R}^{r_p}$ convert the metabolites into macromolecules or vice versa.
\end{itemize}
We are mainly interested into the time evolution of the species and reactions, $X(t) \subset \mathbb{R}^{n_x},~V_y(t) \subset \mathbb{R}^{r_y}$, etc.
But for convenience and easy readability we drop the dependency on time $t$ unless it is important.
The effect of a flux $V_i$ on the time evolution of the species is given by the stoichiometric matrix $S \in \mathbb{R}^{n_y+n_x+n_p,n_r}$ as
\begin{align}\label{eq:original_ode}
\frac{\mathrm{d}}{\mathrm{d}t}\begin{pmatrix}
 Y(t) \\ X(t) \\ P(t)
\end{pmatrix} = \begin{pmatrix}
 \dot{Y}(t) \\ \dot{X}(t) \\ \dot{P}(t)
\end{pmatrix} = S \begin{pmatrix} V_y (t)\\ V_x (t)\\ V_p(t) \end{pmatrix} = \begin{pmatrix} S_y^y V_y(t) \\ S_y^x V_y(t) + S_x^x V_x(t) + S_p^x V_p(t) \\ S_p^p V_p(t) \end{pmatrix},
\end{align}
where the matrix entries $s_{i,j}$ describe the stoichiometry of species $i$ in reaction $V_j$.
The macromolecules $P$ are usually composed of a large number of the small metabolites $X$, which may lead to ill-conditioned dynamics \eqref{eq:original_ode} with large coefficients in the matrix $S_p^x$.
We tackle this by scaling the macromolecules and the according fluxes by a large factor $\alpha \in \mathbb{R}$ as
\begin{gather}
 \Ptil  = \alpha P,~ \tilde{V}_p  = \alpha V_p. \label{eq:scaling}
\end{gather}
Furthermore, we assume the internal metabolites to be in a quasi-steady-state $\dot{X}=0$.
This assumption is valid, because the metabolic reactions $V_x$ are much faster in comparison to the biomass reactions $V_p$ and the external species $Y$ are slow variables due to the large external volume (cf. \cite{waldherr2015}).
This gives us the \textit{boundary layer condition} of the scaled system as
\begin{align}\label{eq:boundary_layer_condition}
 S_y^x V_y(t) + S_x^x V_x(t) + \alpha^{-1} S_p^x {\tilde{V}}_p(t) = 0,
\end{align}
so that $S_p^x$ is normalized with $\alpha$.
For convenience we will write $\tilde{V} = (V_y, V_x, {\tilde{V}}_p)$ for the scaled fluxes.
Let us first discuss the other constraints in the model.
The most important one is the \textit{enzyme capacity constraint}.
Most of the macromolecules $P$ are enzymes whose amounts limit the possible reaction fluxes.
We assume that the first $n_m$ entries in the vector $P$ correspond to these enzymes.
The maximal reaction rates are determined by the reaction specific catalytic constant $k_{\mathrm{cat},\pm j}, j\in \{1,\ldots,n_r\}$.
We differentiate between the forward constant $k_{\mathrm{cat},+j}$ and the backward constant $k_{\mathrm{cat},-j}$, if the reaction is reversible and both directions are catalyzed by the same enzyme.
Since some enzymes are capable of catalyzing multiple reactions, we denote the set of reactions catalyzed by the enzyme $P_i$ as
\begin{align}
 \mathrm{cat}(i) = \{ j \in \mathbb{N}~|~P_i \text{ catalyzes }V_j\}.
\end{align}
The constraint for the enzyme $P_i$ then reads in short
\begin{align}\label{eq:ECC_base_constraint}
 \sum_{j \in \mathrm{cat}(i)} \left| \frac{V_j}{k_{\mathrm{cat},\pm j}} \right| \leq P_i.
\end{align}
To eliminate the absolute value in constraint \eqref{eq:ECC_base_constraint} and transform it into a linear constraint, we have to include all possible sign-combination of the occurring $k_{\mathrm{cat}}$-values.
As an illustration consider $P_1$ catalyzes $V_1$ reversibly and $V_2$ irreversibly.
Then we can write \eqref{eq:ECC_base_constraint} as
\begin{align}
 H_{c,1} V &= \begin{pmatrix}
              k_{\mathrm{cat},+1}^{-1} &  k_{\mathrm{cat},+2}^{-1} & 0 & \ldots & 0 \\
             -k_{\mathrm{cat},-1}^{-1} &  k_{\mathrm{cat},+2}^{-1} & 0 & \ldots & 0
             \end{pmatrix} V \leq \begin{pmatrix} 1 & 0 & \ldots & 0 \\ 1 & 0 & \ldots & 0 \end{pmatrix} P = H_{E,1} P.
\end{align}
The concatenations of the matrices $H_{c,1}$  to $H_{c,n_m}$ and $H_{E,1}$ to $H_{E,n_m}$, are written as $H_c$ and $H_E$, respectively.
We include the scaling of $\Ptil$ by splitting this into the different kinds of fluxes as
\begin{align}
 \alpha H_{c,y} V_y + \alpha H_{c,x} V_x + H_c {\tilde{V}}_p \leq H_E \Ptil.
\end{align}
and use the short notation 
\begin{align}\label{eq:ECC}
 \tilde{H}_c {\tilde{V}} - H_E \Ptil \leq 0.
\end{align}
Usually, any organism needs a certain amount of structural macromolecules to keep working, e.g. the cell wall separating it from the outside.
We express this necessity by enforcing certain percentages of the scaled biomass $b^T \Ptil$ to be made of the structural components, e.g for a structural protein $\Ptil_s$
\begin{align}\label{eq:biomass_composition_base}
 \psi_s b^T \Ptil \leq \Ptil_s,
\end{align}
with $0 < \psi_s < 1$ being the minimal fraction of the total biomass for $\Ptil_s$.
As before, the extension of \eqref{eq:biomass_composition_base} to the network level can be expressed by collecting the individual constraints into
\begin{align}\label{eq:biomass_composition}
 H_B \Ptil \leq 0,
\end{align}
where the rows of $H_B$ correspond to $\psi_s b^T - e_s^T$ for different indices $s$.
We call \eqref{eq:biomass_composition} the \emph{biomass composition constraint}.
Moreover, we consider the \textit{positivity constraint} of all species
\begin{align}\label{eq:positivity_constraint}
 Y \geq 0,~ X\geq 0,~\Ptil \geq 0 
\end{align}
and \textit{biomass-independent flux constraints}
\begin{align}\label{eq:flux_constraints}
 V_{\min} \leq V \leq V_{\max}.
\end{align}

The last step to construct the optimization problem is the introduction of the \textit{objective functional} $J$ as
\begin{align}\label{eq:objective_func}
 J(T, Y_0, \Ptil_0) = \int_0^{T} b^T \Ptil(t)~\mathrm{d}t,
\end{align}
depending on the initial values $Y_0 = Y(0),~\Ptil_0=\Ptil(0)$ and the chosen end-time $T > 0$.
Since $b^T \Ptil$ is proportional to the total biomass, this objective functional represents the accumulation of biomass over time and corresponds to the assumption that the modeled organism optimizes itself to grow as fast as possible.
Finally, we can formulate the full optimal control problem as
\begin{align}\label{eq:deFBA_problem}
\begin{split}
 \underset{Y(\cdot), \Ptil(\cdot),{\tilde{V}}(\cdot)}{\max}~ & \int_0^{T} b^T \Ptil(t) ~ \mathrm{d}t \\
 \mathrm{s.t.}     ~ & \begin{pmatrix} \dot{Y}(t) \\ \dot{\Ptil}(t) \end{pmatrix} = \begin{pmatrix} S_y^y & 0 & 0\\ 0 & 0 & S_p^p \end{pmatrix} \begin{pmatrix} V_y(t) \\ {\tilde{V}}_p(t) \end{pmatrix} \\
		    &  Y(0) = Y_0,~\Ptil(0) = \Ptil_0 \\
		    &  S_y^x V_y(t) + S_x^x V_x(t) + \alpha^{-1} S_p^x {\tilde{V}}_p(t) = 0 \\
		    &  \tilde{H}_c {\tilde{V}}(t) - H_E \Ptil(t) \leq 0\\
		    &  H_B \Ptil(t) \leq 0 \\
		    &  Y(t), \Ptil(t) \geq 0 \\
		    &  V_{\min} \leq {\tilde{V}}(t) \leq V_{\max}
\end{split}
\end{align}
Depending on the initial conditions $Y_0,\Ptil_0$, this linear dynamic optimization problem determines the optimal fluxes $\tilde{V}(t)$ and the corresponding enzyme levels $\Ptil(t)$ at all times $t\in [0,T]$.

The problem with this approach is, that the optimization can lead to artificial results.
For example, the optimal solution may only predict the production of the macromolecules with the best biomass yield.
Whereas, this is an optimal solution from a mathematical standpoint, we know that biological systems usually do not behave this way.
We present this behavior in the following example.

\subsection{Example: Enzymatic growth}\label{sec:example_1}
\begin{figure}[t!]
\begin{minipage}{0.32\textwidth}
 \begin{overpic}[scale=0.55, tics=10]{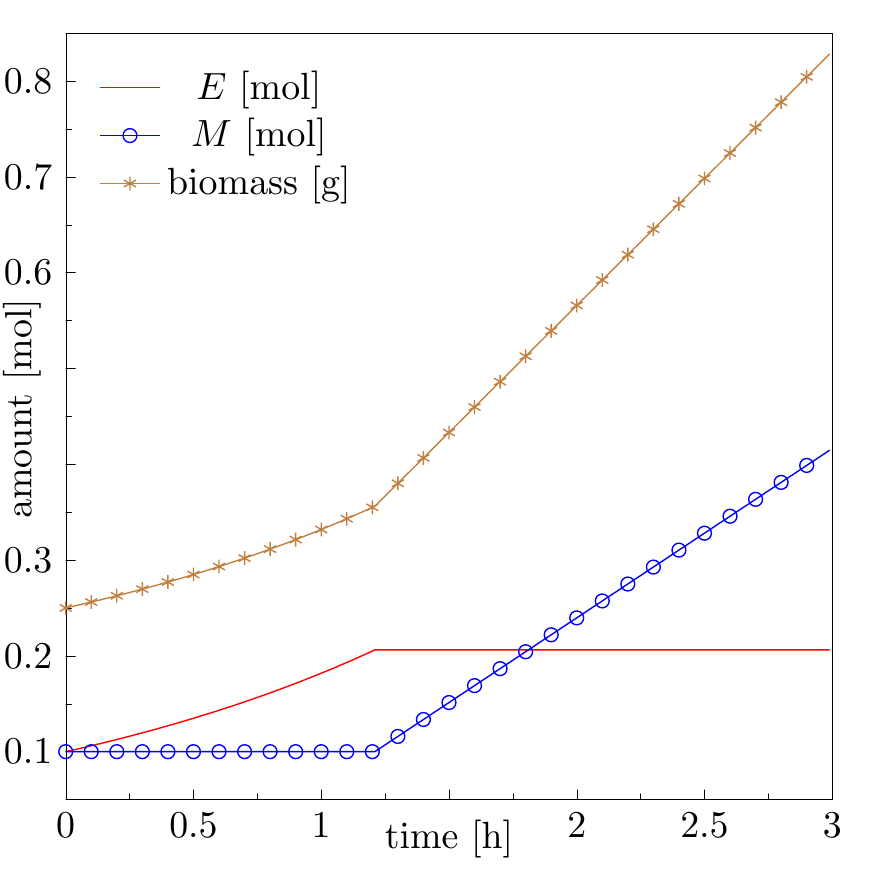}
  \put (50,85) {\textcircled{\smaller A}}
 \end{overpic}
\end{minipage}
\hfill
\begin{minipage}{0.32\textwidth}
 \begin{overpic}[scale=0.55, tics=10]{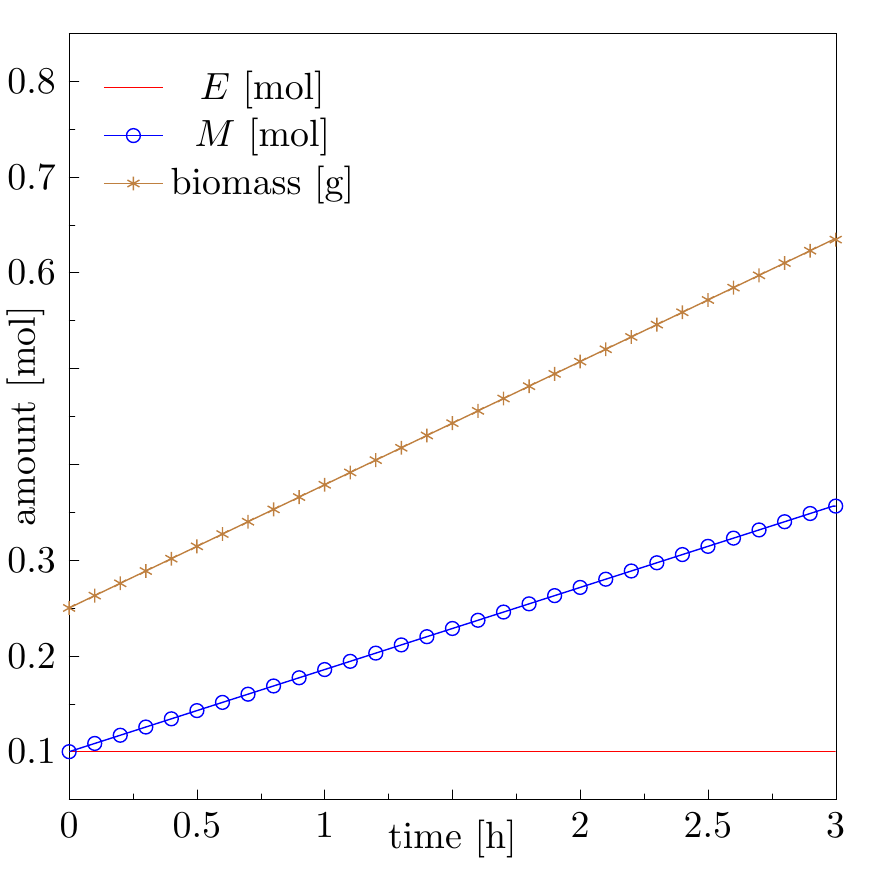}
  \put (50,85) {\textcircled{\smaller B}}
 \end{overpic}
\end{minipage}
\hfill
\begin{minipage}{0.32\textwidth}
 \begin{overpic}[scale=0.55, tics=10]{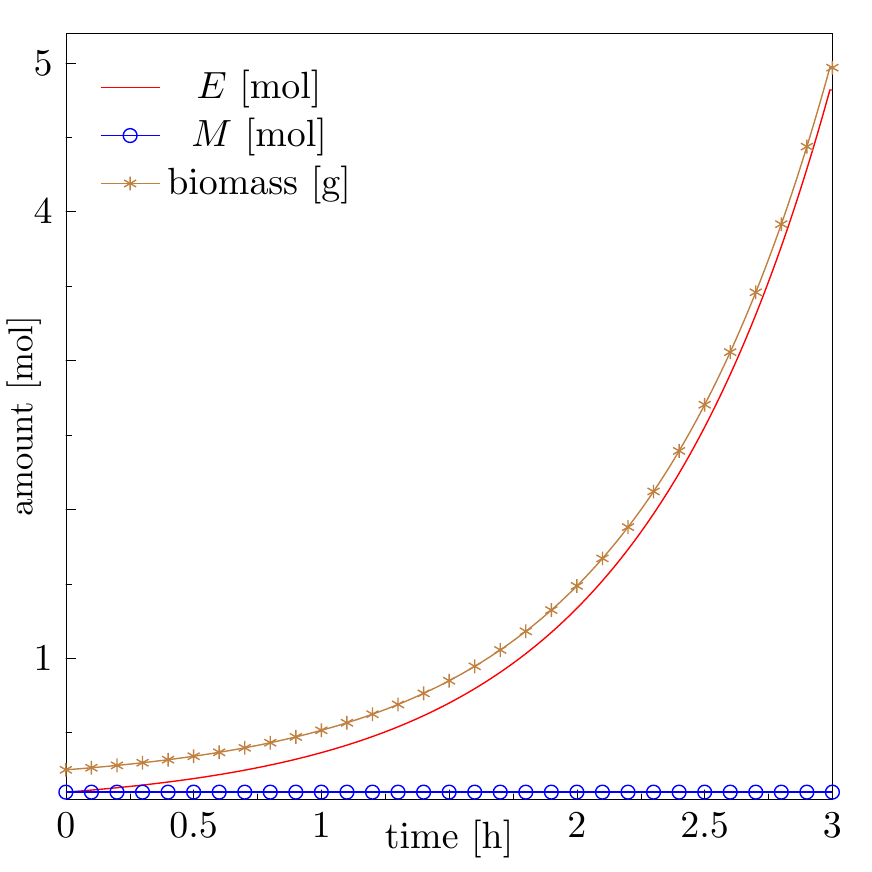}
  \put (50,85) {\textcircled{\smaller C}}
 \end{overpic}
\end{minipage}
\caption{Numerical results for the example in Section \ref{sec:example_1}. Plot {\textcircled{\smaller A}} shows a mixture of linear and exponential growth for $T=3~\mathrm{h}$ and $k_E=1~\mathrm{h}^{-1}$, {\textcircled{\smaller B}} shows pure linear growth for $k_E = 0.5~\mathrm{h}^{-1}$, and {\textcircled{\smaller C}} shows exponential growth until the end-time $T=3$ h for $k_E = 10~\mathrm{h}^{-1}$.}\label{fig:example_1_results}
\end{figure}
We present a simple model, which we already analyzed very detailed in \cite{Waldherr2017}, to give the reader an idea how changing the catalytic constant $k_{\mathrm{cat}}$ or other system parameters, e.g. the end-time $T$, can impact the quality of results.
In this model the organism can either invest in the production of catalyzing enzymes to increase its capability for growth or it can produce storage components, which have a better biomass yield.
We use a minimal model size by introduction of one external nutrient $N$, one internal metabolite $A$, two macromolecules $M,E$, and three irreversible reactions $V_M,~V_A,~ V_E$ with stoichiometry as shown in \eqref{eq:example_1_network1}-\eqref{eq:example_1_network3}.
\begin{subequations}
\begin{eqnarray}
 V_A : & N            & \rightarrow  A \label{eq:example_1_network1} \\ 
 V_E : & 100~ N + 100 ~ A & \rightarrow 1 E \label{eq:example_1_network2} \\
 V_M : & 100~ N + 100 ~ A & \rightarrow 1 M. \label{eq:example_1_network3}
\end{eqnarray}
\end{subequations}
The external nutrient $N$ represents a collection of components necessary for growth, such as carbon, nitrogen, etc.
Further processed components made from these nutrients are collected as the internal metabolite $A$.
We differentiate the macromolecules into the group of enzymes $E$, collecting the whole enzymatic machinery needed for growth, and non-enzymatic macromolecules $M$.
These reflect storage components such as lipids, starch, and glycogen.
We handle the large difference in the stoichiometric values by scaling with the factor $\alpha=100$ as described in \eqref{eq:scaling}.
The scaled dynamics for this model are
\begin{align}
\frac{\mathrm{d}\tilde{x}}{\mathrm{d}t}  = \tilde{S} {\tilde{V}} \Leftrightarrow \frac{\mathrm{d}}{\mathrm{d}t}
 \begin{pmatrix}
  x_N \\ x_A \\ \tilde{x}_E \\ \tilde{x}_M
 \end{pmatrix} & = \begin{pmatrix} -1 & -1 & -1  \\
				    1 & -1 & -1 \\
				    0 & 1 & 0 \\ 
				    0 & 0 & 1 \end{pmatrix} \begin{pmatrix} V_A \\ {\tilde{V}}_E \\ {\tilde{V}}_M \end{pmatrix} 
\end{align}
According to steady-state assumption \eqref{eq:boundary_layer_condition}, we do not allow accumulation of the metabolites $A$
\begin{align}
 & \dot{x}_A  = 0 =  V_A - {\tilde{V}}_E - {\tilde{V}}_M.
\end{align}
We further assume all three reactions are catalyzed by the enzymatic macromolecules $E$ with the respective catalytic constants $k=(k_E,~k_M,~k_A)$.
Hence, the enzyme capacity constraint \eqref{eq:ECC} for the scaled system is given as
\begin{align}
 \alpha \frac{V_A}{k_A}+ \frac{{\tilde{V}}_E}{k_E} + \frac{{\tilde{V}}_M}{k_M}  \leq \tilde{x}_E.
\end{align}
We collect the molecular weights of the species in the biomass vector $b^T = (b_E, b_M)$ and define the objective functional as
\begin{align}
 \max_{{\tilde{V}}(\cdot)} J({\tilde{V}}(\cdot), \tilde{x}_0) = \max_{{\tilde{V}}(\cdot)} \int_0^T b^T \tilde{x}(t) ~\mathrm{d} t =  \max_{{\tilde{V}}(\cdot)}\int_0^T b_M \tilde{x}_M(t) + b_E \tilde{x}_E(t) ~\mathrm{d} t,
\end{align}
with initial values $\tilde{x}_0$ and fixed end-time $T$. 
We have already shown analytically in \cite{Waldherr2017} that this system can experience basically three different growth modes depending on the chosen parameters and end-time $T$.
Growth in this sense means the change of the biomass over time $b^T \dot{x}$.
\begin{itemize}
 \item A \emph{stationary phase} with $\dot{x} = 0$, which occurs as soon as nutrients deplete.
 \item A \emph{linear growth} phase starting at $\tau$ with the optimal solution 
 \begin{align}
   \begin{split}
  x^*_N (t) &= x^*_N(\tau) - 2 \tilde{x}^*_M (t) \\
  \tilde{x}^*_E (t) &= \tilde{x}_E (\tau)\\
  \tilde{x}^*_M (t) &= \tilde{x}_M (\tau) + \frac{k_A k_M}{k_A + \alpha k_M} \tilde{x}_E(\tau) (t-\tau),
   \end{split}
 \end{align}
 meaning the solution predicts only increase in the storage $M$, while keeping the amount of enzymes $E$ constant.
 \item An \emph{exponential growth} phase starting at $\tau$ with the optimal solution
 \begin{align}
   \begin{split}
    x^*_N (t) & = x^*_N(\tau)  - 2 x^*_E (t) \\
    \tilde{x}^*_E (t) & = \tilde{x}_E(\tau) \exp\left( \frac{k_A k_E}{k_A + \alpha k_E} (t - \tau) \right)\\
    \tilde{x}^*_M (t) & = \tilde{x}_M(\tau).
   \end{split}
 \end{align}
\end{itemize}
If we assume unlimited nutrients, we would expect a biological system to work exclusively in the exponential phase and produce no storage $M$ at all.
But the formulation as an optimization problem can generate solutions containing linear growth phases depending on the system parameters and the end-time.
For this example, we fix the end-time $T=3$ h, and vary the value of $k_E \in [0.5\text{ h}^{-1}, 5\text{ h}^{-1}]$ with a nominal value of $k_E = 1 \text{ h}^{-1}$.
The other parameters are shown in Table \ref{tb:example_1}.
The initial amounts are chosen as $x_N(0)=10^6$ mol, $x_E(0) = 0.1$ mol, and $x_M(0) = 0.1$ mol to ensure nutrients are not limiting.
We refrain from plotting the results for $N$ as these are given by the dynamics
\begin{align}
 x^*_N (t) & = x^*_N(0)  - 2 ( \tilde{x}^*_E (t) + \tilde{x}^*_M(t) ).
\end{align}

We can observe a partwise linear solution plotted in Figure \ref{fig:example_1_results} {\textcircled{\smaller A}}.
As mentioned, we chose the parameters such that the yield of the storage $M$ is larger ($b_M > b_E$), but depending on the parameters exponential growth is preferred as it enables quicker nutrient uptake.
Up until $t=1.12$ h, the solution shows an exponential increase in enzymes, but afterwards all resources are used to maximize biomass by producing only storage.
This contradicts our expectation of exclusive exponential growth in a non-limiting nutrient situation, but slight changes to $k_E$ can change these results completely.
Using the minimal value $k_E = 0.5~\mathrm{h}^{-1}$, we can observe a solution consisting of a single linear growth phase in {\textcircled{\smaller B}}.
This change in the reproductive capabilities of $E$ is enough to make purely linear growth the most effective solution.
Following \cite{Waldherr2017}, we can even pinpoint the exact conditions necessary for optimal linear growth
\begin{align}
 T &< \frac{2(k_M b_M - k_E b_E)}{b_M k_M k_E} \text{ and } x_N(T) \geq 0 \label{eq:condition_linear_growth_example_1}
\end{align}
Equation \eqref{eq:condition_linear_growth_example_1} also shows, that the end-time $T$ is as important as the other system parameters for the quality of the solution.
The right plot {\textcircled{\smaller C}} uses the maximal value of $k_E = 10~\mathrm{h}^{-1}$ and shows a purely exponential growth behavior.
This is interesting because the instantaneous yield of the storage $M$ is still better in this case ($b_M > b_E$), however, the exponential solution is able to outgrow the linear growth mode for large $k_E$ values.

Inspecting the results of this minimal example, we are now facing two problems. 
Can we eliminate the mathematical artifacts from solutions like the one presented in {\textcircled{\smaller A}} and how can we handle the uncertain value $k_E$.
Interestingly, the answers to these questions are highly related.
First we extend the deFBA with means from Model Predictive Control (MPC) to answer the first one.

\begin{table}
\caption{Numerical values for the parameters in Example 1}\label{tb:example_1}
\centering
\begin{tabular}{lccccccc}
\toprule
parameter & $b_M~[\mathrm{g}~\mathrm{mol}^{-1}]$ & $b_E~[\mathrm{g}~\mathrm{mol}^{-1}]$ & $k_A~[\mathrm{h}^{-1}]$ & $k_M~[\mathrm{h}^{-1}]$ & $k_E~[\mathrm{h}^{-1}]$ \\
\midrule
value &  150  & 100 & 150 & 2 & 1  \\
\bottomrule
\end{tabular}
\end{table}

\section{Short term deFBA}\label{sec:stdeFBA}
The previous example shows artifacts of the optimization technique in form of linear growth phases (cf. Figure \ref{fig:example_1_results} {\textcircled{\smaller A}} and {\textcircled{\smaller B}}).
We present in this section an extension of the deFBA method, which helps to eliminate these artifacts and can reduce the computational cost of finding the optimal solution.

\subsection{Using a receding prediction horizon}
An obvious mismatch between our models and a real microorganism is a vast difference in knowledge on the development of the environment.
The deFBA depicts a deterministic environment and solves an optimization problem with exact information for all times $t \in [0,T]$.
Real biological systems on the other hand are acting more like Model Predictive Controllers (MPC) \cite{rawlings2009}, meaning they react and adapt to changing conditions on the fly.
Hence, we wanted to use the basic idea of a reacting system in our model by the inclusion of the receding prediction horizon.

In the original deFBA formulation \eqref{eq:deFBA_problem} we solve one large optimization problem over until the end-time $T$.
here, we define instead a \textit{prediction horizon} $p > 0$ representing the time the system can plan ahead.
Thus, we solve an iterative scheme opposed to a single optimization problem.
We introduce a \emph{step size} $h > 0 ,~ h \leq p$ to discretize time on the equidistant grid $\Delta t$
\begin{align}
 \Delta t = \{ t_k = k h~\vert~k=0,\ldots,N\} \label{eq:time_grid}.
\end{align}
Each of the individual optimization problems over the \textit{control times} $kh$, $k \in \{0,\ldots, N\},~kN=T$ is defined as
\begin{align}\label{eq:sdeFBA_continous}
\begin{split}
\underset{Y(\cdot),\Ptil(\cdot),{\tilde{V}}(\cdot)}{\max} &\int_{kh}^{kh+p} b^T \Ptil(t)~\mathrm{d}t \\
 \text{s.t.} &  \begin{pmatrix} \dot{Y}(t) \\ \dot{\Ptil}(t) \end{pmatrix} = \begin{pmatrix} S_y^y & 0 & 0\\ 0 & 0 &S_p^p \end{pmatrix} {\tilde{V}}(t) \\
	     &  Y(kh) = Y_k,~\Ptil(kh) = \Ptil_k \\
	     &  S_y^x V_y(t) + S_x^x V_x(t) + \alpha^{-1} S_p^x {\tilde{V}}_p(t) = 0 \\
	     &  \tilde{H}_c {\tilde{V}}(t) - H_E \Ptil(t) \leq 0\\
	     &  H_B \Ptil(t) \leq 0 \\
	     &  Y(t), \Ptil(t) \geq 0 \\
	     &  V_{\min} \leq {\tilde{V}}(t) \leq V_{\max}. \\
\end{split}
\end{align}
After solving the problem \eqref{eq:sdeFBA_continous} for fixed $k$, the first part $t \in (kh , kh+p)$ of the calculated trajectories is appended to the solution curves $Y^*, \Ptil^*, {\tilde{V}}^*$.
Then the next iteration is initialized with the updated values $Y_{k+1}=Y((k+1)h), \Ptil_{k+1}=\Ptil((k+1)h)$.
For our numerical results we discretize all variables on $\Delta t$ \eqref{eq:time_grid} and denote the approximation of the species states $Y_k, \Ptil_k$ and reaction fluxes ${\tilde{V}}_k$.
We further express the prediction horizon as $p=\hat{p}h$.
The dynamics of the system are substituted by a collocation method $f: \mathbb{R}^{n_y} \times \mathbb{R}^{n_p} \times \mathbb{R}^{n_r} \rightarrow \mathbb{R}^{n_y + n_p}$. 
As an example, the problem \eqref{eq:sdeFBA_continous} can be discretized with an 1-step collocation method $f$ as
\begin{align}\label{eq:sdeFBA}
  \begin{split}
 \underset{(Y_l , \Ptil_l,{\tilde{V}}_{l})_{l \in \{k , \ldots, k + \hat{p} - 1\}}}{\max}~ & \sum_{l = k}^{k+\hat{p}} b^T \Ptil_l\\
 \text{with given}     ~& Y_k, \Ptil_k \\
 \text{and}~\forall l \in\{ k,\ldots, k+\hat{p}-1 \} & \\
		    & \begin{pmatrix} Y_{l+1} \\ \Ptil_{l+1}  \end{pmatrix} = f(Y_{l}, \Ptil_{l}, {\tilde{V}}_{l}) \\
		    &  \begin{pmatrix} S_y^x,~ S_x^x,~ \alpha^{-1} S_p^x \end{pmatrix} {\tilde{V}}_{l} = 0 \\
		    &  \tilde{H}_c {\tilde{V}}_l - H_E \Ptil_l \leq 0\\
		    &  H_B \Ptil_{l+1} \leq 0 \\
		    &  Y_{l+1}, \Ptil_{l+1} \geq 0 \\
		    &  V_{\min} \leq {\tilde{V}}_l \leq V_{\max}
\end{split}
\end{align}
We call the problem \eqref{eq:sdeFBA} the \textit{short-term deFBA} (sdeFBA).
Solving the system starting at the $k$-th time step, gives us the optimal fluxes for the next $\hat{p}$-time steps.
However, we only implement the first fluxes $\{ {\tilde{V}}_{k} \}$ and reformulate \eqref{eq:sdeFBA} for the next iteration for $k+1$.

This approach is beneficial in multiple ways, e.g. we already mentioned the possibility to eliminate mathematical artifacts, see also Example \ref{sec:example_3}.
But the most relevant advantage for us is the availability of algorithms to handle uncertainties in such problems via multi-stage NMPC \cite{lucia2013}.

Since we have already shortly discussed in Example \ref{sec:example_1} that different end-times $T$ may lead to qualitatively different results, we are facing the same problem in choosing the prediction horizon $p$.
In most applications of MPC the prediction horizon is dictated by real-time requirements for the computation and is chosen as large as possible to maximize performance.
For our approach we are interested in a prediction horizon $p_{\mathrm{up}}$ such that a solution including exponential growth phases is attainable.
This means, we need to choose $p_{\mathrm{up}}$ sufficiently large such that linear solutions become suboptimal, whilst keeping $p_{\mathrm{up}}$ as small as possible for the sake of low computational cost.

\subsection{Choosing the prediction horizon}\label{sec:upper_bound}
\begin{figure}
\centering
\includegraphics[scale=0.9]{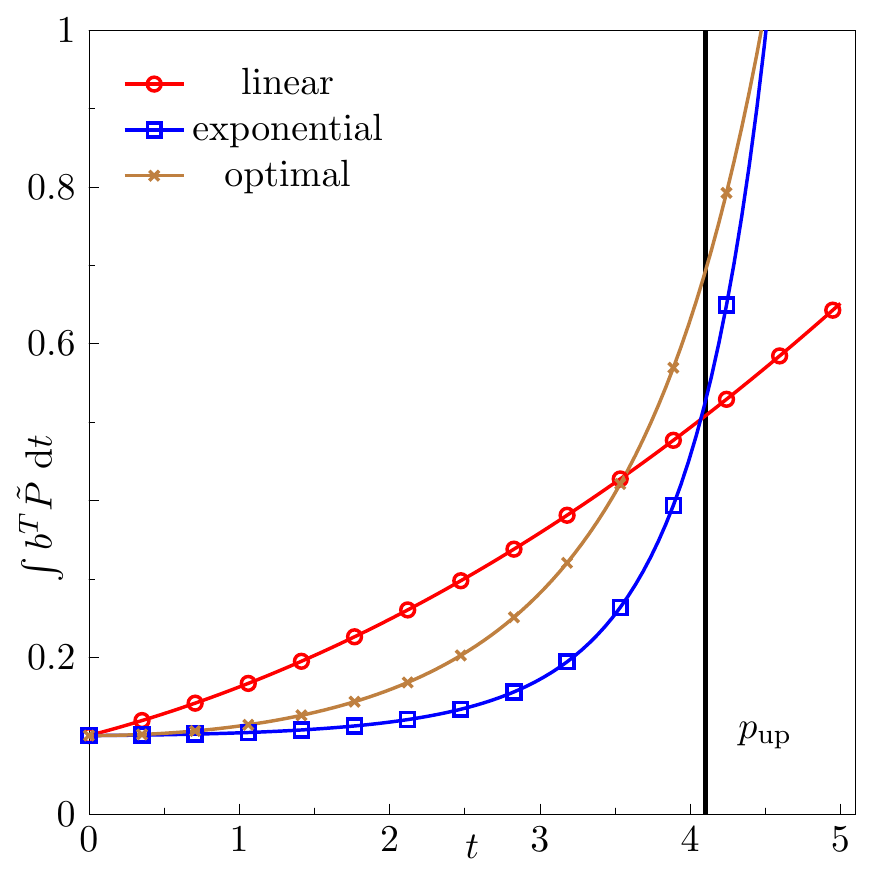}
\caption{Illustration of method to choose the prediction horizon. Upper bound on linear growth shown in red ($\circ$), lower bound on exponential growth in blue ($\square$), and the optimal solution in brown (x).}\label{fig:exp_vs_lin}
\end{figure}
Our idea to choose the prediction horizon is based on the results of Example \ref{sec:example_1}.
Looking at the plots {\textcircled{\smaller A}} and {\textcircled{\smaller B}} from Figure \ref{fig:example_1_results}, we can see that for short times linear growth is faster than exponential growth.
But an exponential solution outgrows a linear one after enough time, as sketched in Figure \ref{fig:exp_vs_lin}. 
Hence, we calculate the maximal achievable linear growth rate and compare it to a suboptimal solution with enforced exponential growth.
The time at which the suboptimal exponential solution outgrows the best linear one is then chosen as the prediction horizon, to ensure the solver will choose an exponential growth solution if available.
In the following we show how to compute these solutions.

Based on the initial biomass composition $\Ptil_0=\Ptil(0)$ we define an upper bound on linear growth $b^T\Ptil_{\mathrm{lin}}(t)$ by maximizing the slope of the curve $b^T \Ptil_{\mathrm{lin}}(t)$ at $t=0$.
Additionally, we set the system dynamics to $\dot{\Ptil} = S^p_p {\tilde{V}}_{\mathrm{lin}}$, with ${\tilde{V}}_{\mathrm{lin}}(\Ptil_0) \in \mathbb{R}^{n_r}$.
This translates to the optimization problem similar to the standard FBA \cite{varma1994metabolic}
\begin{align}
 \begin{split}\label{eq:upper_lin_bound}
  {\tilde{V}}_{\max} = \arg\underset{{\tilde{V}}_{\mathrm{lin}}}{\max}~~ & b^T S_p^p {\tilde{V}}_{\mathrm{lin},p}\\
  \text{s.t.~~} & \begin{pmatrix} S_y^x,~ S_x^x,~ \alpha^{-1} S_p^x \end{pmatrix} {\tilde{V}}_{\mathrm{lin}} = 0 \\
	      & \tilde{H}_c {\tilde{V}}_{\mathrm{lin}} - H_E \Ptil_0 \leq 0 \\
	      & V_{\min} \leq {\tilde{V}}_{\mathrm{lin}} \leq V_{\max}.
 \end{split}
\end{align}
Please note that \eqref{eq:upper_lin_bound} is independent of the real nutrient situation $Y(t)$ as it represents a strict upper bound.
We use the solution ${\tilde{V}}_{\max}$ to obtain the corresponding biomass trajectory by solving the differential equation
\begin{align}\label{eq:lin_grow_ode}
 \dot{\Ptil}_{\mathrm{lin}} &= S_p {\tilde{V}}_{\max} ~ \Rightarrow ~ \Ptil_{\mathrm{lin}}(t) = \Ptil_0 + t S_p {\tilde{V}}_{\max}.
\end{align}
%


To calculate a lower bound on the exponential growth solution we construct another optimization problem inspired by the Resource Balance Analysis (RBA) \cite{goelzer2011cell}.
In RBA one is looking for enzyme levels which maximize the growth rate $\mu(\Ptil_0) \in \mathbb{R}$ of the system.
We enforce balanced growth \cite{baroukh2014drum}, meaning the biomass composition is fixed during growth, via the dynamics $\dot{\Ptil} = \mu \Ptil$.
As before, we solve this only at time $t=0$ and identify the flux vector ${\tilde{V}}_{\exp} \in \mathbb{R}^{n_r}$ realizing this growth rate.
\begin{align}
 \begin{split}\label{eq:lower_exp_bound}
  \mu_{\max} = \underset{{\tilde{V}}_{\exp},\mu}{\max}~~ & \mu\\
  \text{s.t.~~}& \mu \Ptil_0 =(0,0,S_p^p) {\tilde{V}}_{\exp} \\
	      & \begin{pmatrix} S_y^x,~ S_x^x,~ \alpha^{-1} S_p^x \end{pmatrix} {\tilde{V}}_{\mathrm{exp}} = 0 \\ 
              & \tilde{H}_c {\tilde{V}}_{\exp} - H_E \Ptil_0 \leq 0 \\
	      &  V_{\min} \leq {\tilde{V}}_{\exp} \leq V_{\max} \\
	      & \Ptil \geq 0.
 \end{split}
\end{align}
The biomass curve with the maximal growth rate $\mu_{\max}$ can be identified by solving the system dynamics
\begin{align}\label{eq:exp_grow_ode}
 \dot{\Ptil} &= \mu_{\max} \Ptil ~ \Rightarrow ~ \Ptil(t) = \Ptil_0 \exp(\mu_{\max} t).
\end{align}
With solutions for linear and exponential growth in hand, we now calculate the contact time $p_{\mathrm{up}}$, at which the integrated biomass curves meet.
The integral of the biomass curve for exponential growth is derived from \eqref{eq:exp_grow_ode} as
\begin{align}
 B_{\exp}(\mu_{\max}, p_{\mathrm{up}}, \Ptil_0) = \int_0^{p_{\mathrm{up}}} b^T \Ptil(t) ~\mathrm{d}t &= \mu_{\max}^{-1} b^T \Ptil_0 \exp(\mu_{\max} {p_{\mathrm{up}}}) - \mu_{\max}^{-1} b^T \Ptil_0 \label{eq:exp_solution}
\end{align}
and the corresponding integral for the linear case \eqref{eq:lin_grow_ode} is computed as
\begin{align}
 B_{\mathrm{lin}}({p_{\mathrm{up}}},{\tilde{V}}_{\max},\Ptil_0) = \int_0^{p_{\mathrm{up}}} b^T \Ptil_{\mathrm{lin}}(t) ~\mathrm{d}t &= {p_{\mathrm{up}}} b^T \Ptil_0  + {p_{\mathrm{up}}}^2 b^T S {\tilde{V}}_{\max}. \label{eq:linear_solution}
\end{align}
Both solutions \eqref{eq:exp_solution} and \eqref{eq:linear_solution} satisfy the deFBA constraints \eqref{eq:deFBA_problem} but usually are suboptimal solutions to the original deFBA problem.
Therefore, we interpret \eqref{eq:exp_solution} as a lower bound for exponential growth.
We calculate the contact time $p_{\mathrm{up}}$ by solving 
\begin{align}\label{eq:determine_pexp}
 B_{\mathrm{lin}}(p_{\mathrm{up}}) = p_{\mathrm{up}} b^T \Ptil_0  + p_{\mathrm{up}}^2 b^T (0,0,S_p^p) {\tilde{V}}_{max} = \mu^{-1} b^T \Ptil_0 \exp(\mu p_{\mathrm{up}}) - \mu^{-1} b^T \Ptil_0 = B_{\exp}(p_{\max}).
\end{align}
The solution $p_{\mathrm{up}}=0$ always exists for \eqref{eq:determine_pexp}.
If this is the only solution, the sub-optimal biomass curve $B_{\exp}$ is larger at all times and we can choose the prediction horizon arbitrarily and still get an exponential solution from the sdeFBA.
Secondly, there may exists another solution for \eqref{eq:determine_pexp} with $p_{\mathrm{up}} > 0$.
This implies that on a horizon $p \geq p_{\mathrm{up}}$ the sub-optimal exponential solution outgrows the best linear solution (cf. Figure \ref{fig:exp_vs_lin}).
Thus, as long as nutrients do not deplete and we choose $p \geq p_{\mathrm{up}}$, the optimal solution for \eqref{eq:sdeFBA_continous} must be an exponential growth phase.
Because we calculated $p_{\mathrm{up}}$ based on a comparison of upper and lower bounds without the inclusion of the actual nutrient dynamics, depending on the specific situation it could be that a smaller prediction horizon also gives exponential growth.

\section{Robust deFBA}\label{sec:rdeFBA}
\subsection{Using multi-stage MPC}
Having the sdeFBA established in the last section, we can now show how we include uncertainties in our method via the multi-stage NMPC approach presented in \cite{lucia2013, lucia2014}; also presented in e.g. \cite{delaPena2005b} for linear systems.
As mentioned before we want to model very different uncertainties, but mostly unknowns in the model structure like parameter uncertainties, unknown stoichiometry, etc. 
While our framework can handle different kinds of uncertainties at the same time, we focus in this work on the inclusion of measurement errors in the $k_{\mathrm{cat}}$-values as leading example.

\begin{figure}
\centering
\includegraphics[scale=.95]{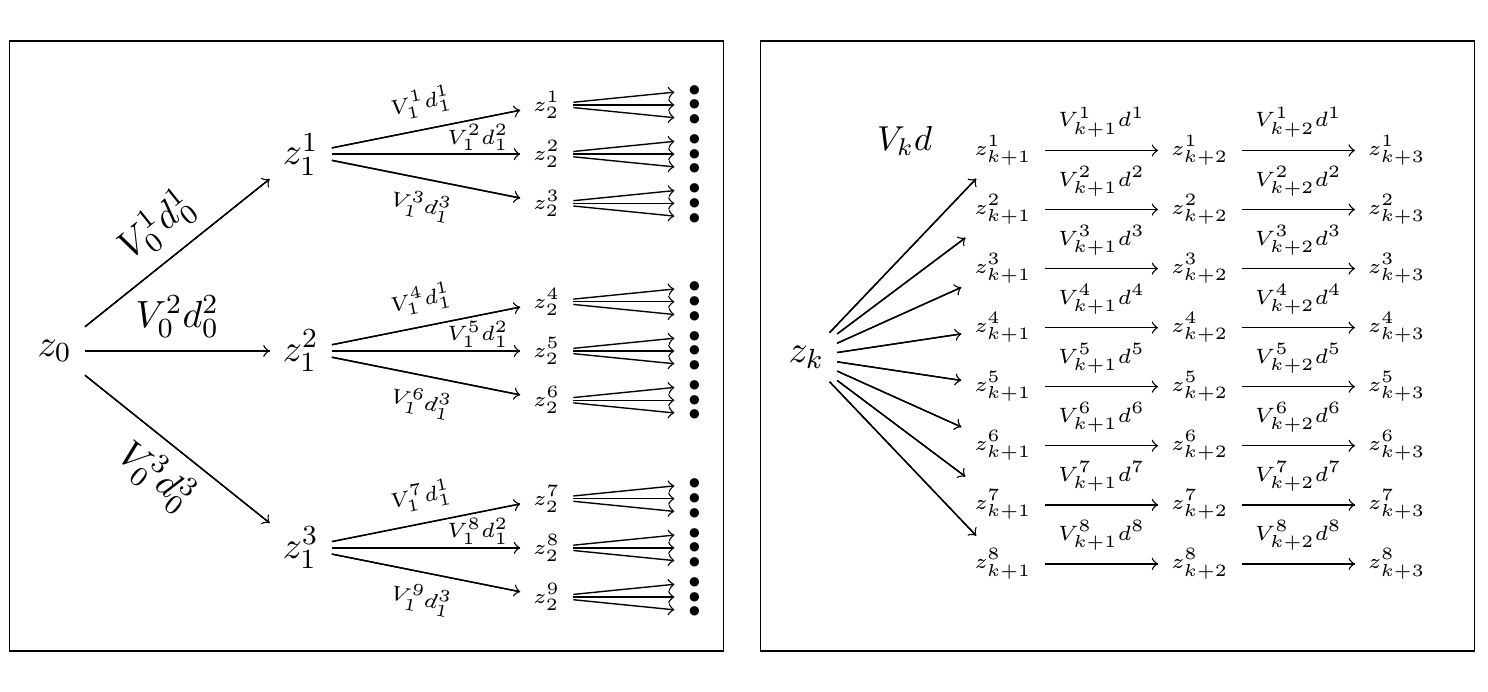}
\caption{Scenario trees for different robust horizons. (Left) Example of a full branching tree. (Right) Reduced scenario tree for a robust horizon $r=1$ and three uncertainties. Only maximal deviations from nominal values are included.}\label{fig:scenario_tree}
\end{figure}

The multi-stage NMPC is based on a representation of the uncertainties by a branching tree as shown in Figure \ref{fig:scenario_tree} (left).
Each branch represents the effect of an uncertainty $d^j$ and the optimal fluxes $V^j$ to handle this effect.
As we can only handle a finite amount of scenarios we assume all uncertainties are bounded $d^j(t) \in [d^j_{\min},d^j_{\max}]$ at all times.
A path in the tree from the root to a final leaf is called a \emph{scenario}.
To make the notation easier to read, we use an upper index for scenarios and a lower index for prediction steps.
So with each new iterative step $k$, we get another set of scenarios depending on the number of uncertainties $n_d$.
The main challenge of this approach is that the scenario tree grows exponentially with the prediction horizon $p$.
However, previous results \cite{lucia2014} have shown that a full coverage of the scenario tree is usually not necessary and we can assume that the uncertainties stay constant after a certain time span.
We choose this time span to be identical to the step size $h$ of the time grid used for the discrete approximation of the system.
This way we can simplify the tree to the form presented in Figure \ref{fig:scenario_tree} (right).
Constructing the scenarios with the extreme values $d^j_{\min},d^j_{\max}$ for each $d^j$, we can ensure the calculated fluxes $V_0$ are feasible for all scenarios.
Thus, considering $n_d$ uncertainties we get a total of $2^{n_d}$ different scenarios, with the length of each scenario being defined by the prediction horizon $p$.
Figure \ref{fig:scenario_tree} (right) shows how to handle a system with three uncertainties where each branch corresponds to a $d^j$ of the set
\begin{align}\begin{split}
 d^j \in \Big\{&  (d^1_{\min},d^2_{\min},d^3_{\min}), (d^1_{\min},d^2_{\min},d^3_{\max}), (d^1_{\min},d^2_{\max},d^3_{\min}),(d^1_{\min},d^2_{\max},d^3_{\max}),\\
	       &(d^1_{\max},d^2_{\min},d^3_{\min}),(d^1_{\max},d^2_{\min},d^3_{\max}),(d^1_{\max},d^2_{\max},d^3_{\min}),(d^1_{\max},d^2_{\max},d^3_{\max}) \Big\}.
	       \end{split}
\end{align}
The chronological first set of fluxes $V_0$ in each iteration is the most important one, as this is the one we want to implement to the final solution in each iteration.
These fluxes have to be feasible under all uncertainties $d^j$, which we denote with $d$.

Multi-stage MPC dictates us to implement the scenario tree into the discrete sdeFBA \eqref{eq:sdeFBA} by optimizing over all scenarios at once.
The new objective function is chosen as the weighted sum of the individual objectives \eqref{eq:objective_func}
\begin{align}
 \hat{J}(t_0,t_f) = \sum\limits_{j=1}^{2^{n_d}} \omega^j J^j = \sum\limits_{j=1}^{2^{n_d}}\int_{t_0}^{t^f} \omega^j b^T \Ptil^j~\mathrm{d}t, \label{eq:robust_objective_func}
\end{align}
with start-time $t_0$, end-time $t_f$, and the weights $\omega^j$.
The weights can be used if additional information on the likelihood of a scenario is given, but usually we use $\omega^j = 1$.
Further, we define an index set ${\cal{J}} = \{1,2, \ldots,  2^{n_d} \}$ for the sake of clarity.
The robust approach is then defined as an iterative scheme over $k \in \{0, \ldots , N\}$
\begin{subequations}\label{eq:rdeFBA}
\begin{align}\label{eq:rdeFBA_part1}
 \underset{ (Y^j, \Ptil^j, {\tilde{V}}^j)_j, j \in {\cal{J}}}{\max}~ & \sum_{j=1}^{2^{n_d}}\int\limits_{kh}^{kh+p} b^T \Ptil^j~\mathrm{d}t\\
 \text{subject to }\forall j \in {\cal{J}}:~ & \begin{pmatrix} \dot{Y}^j \\ \dot{\Ptil}^j \end{pmatrix} = \begin{pmatrix} S_y^y & 0 & 0\\ 0 & 0 & S_p^p \end{pmatrix} {\tilde{V}}^j \\
		    & Y(hk) = Y_k ,\Ptil(hk)=\Ptil_k & \\
		    & \begin{pmatrix} S_y^x,~ S_x^x,~ \alpha^{-1} S_p^x \end{pmatrix} {\tilde{V}}^{j} = 0 \\
		    &  \tilde{H}^j_c {\tilde{V}}^j - H_E \Ptil^j \leq 0\\
		    &  H_B \Ptil^j \leq 0 \\
		    &  Y^j, \Ptil^j \geq 0 \\
		    &  V_{\min} \leq {\tilde{V}}^j \leq V_{\max} \label{eq:rdeFBA_part8}
\intertext{Additionally we have to include the \emph{non-anticipativity constraint}, which enforces the copies of the first part of the trajectories to be identical for all scenarios.}
 \forall~a,b \in {\cal J} &,~ {\tilde{V}}^a(t) = {\tilde{V}}^b(t) = {\tilde{V}}(t),~\forall t \in (kh, (k+1)h)  \label{eq:non_antip_con}.
\end{align}
\end{subequations}
The first part of the problem \eqref{eq:rdeFBA_part1}-\eqref{eq:rdeFBA_part8} is identical to the approach we presented for the sdeFBA, but optimizes over all scenarios at once.
We call the full scheme \eqref{eq:rdeFBA} the \emph{robust deFBA} (rdeFBA).
In the case of measurement errors, the only influence of the uncertainties is found in the matrices $\tilde{H}^j_c$, but the framework allows to make all entries in the constraints dependent on the scenario.
As with the sdeFBA, we implement only the first time interval ${\tilde{V}}(t),t \in (kh, (k+1)h)$ to the solution.
For both examples presented in this work, the dynamics are discretized using orthogonal collocation on finite elements using the \textsf{do-mpc} tool \cite{lucia2016_CEP}.
Within \textsf{do-mpc} the resulting optimization problems are solved by IPOPT \cite{wachter2006}, with exact first and second order derivatives computed via CasADi \cite{andersson2012}.
\textsf{Do-mpc} also offers an efficient implementation of multi-stage NMPC.
While the rdeFBA constructs linear programs, the tools we are using support non-linear problems and we can easily include non-linear dynamics or constraints in our models.

We calculate the prediction horizon for the robust deFBA in the way suggested in Section \ref{sec:upper_bound}.
Keep in mind that changes in the catalytic constants may lower the achievable growth rate of the system significantly.
Hence, it may be necessary to calculate the prediction horizon in the scenario $d^1$ with minimal $k_{\mathrm{cat}}$-values.

\subsection{Example: Enzymatic Growth}\label{sec:example_3}
\begin{figure}[t!]
\begin{minipage}{0.45\textwidth}
 \begin{overpic}[scale=0.75, tics=10]{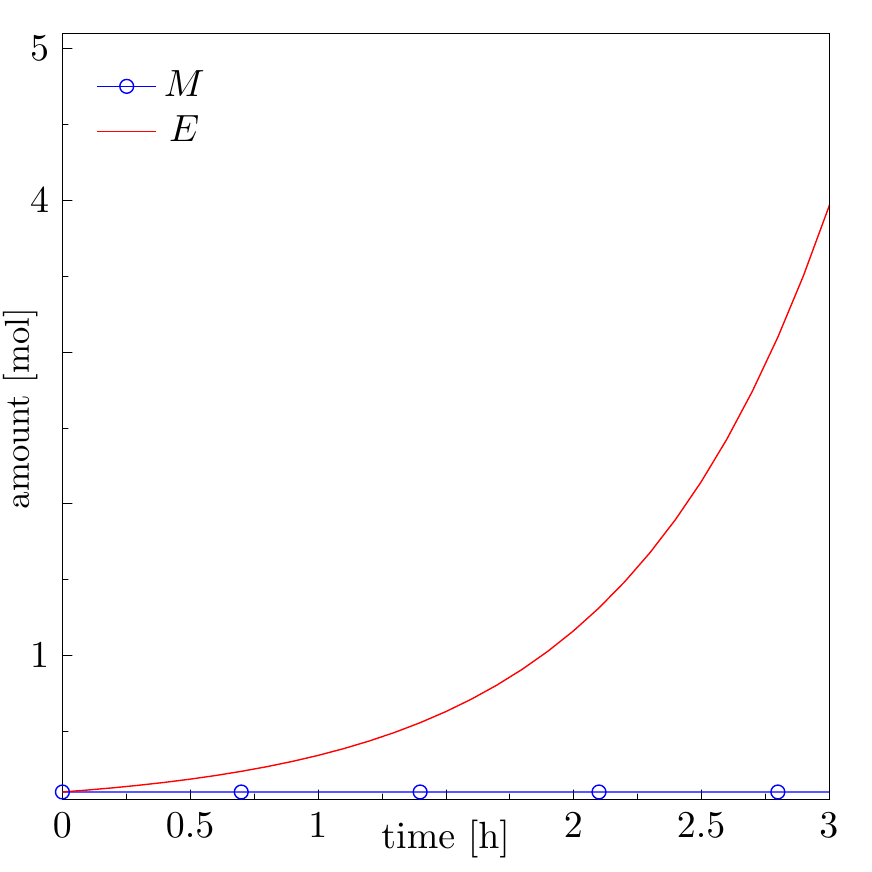}
  \put (50,85) {\textcircled{\smaller A}}
 \end{overpic}
\end{minipage}
\hfill
\begin{minipage}{0.45\textwidth}
 \begin{overpic}[scale=0.75, tics=10]{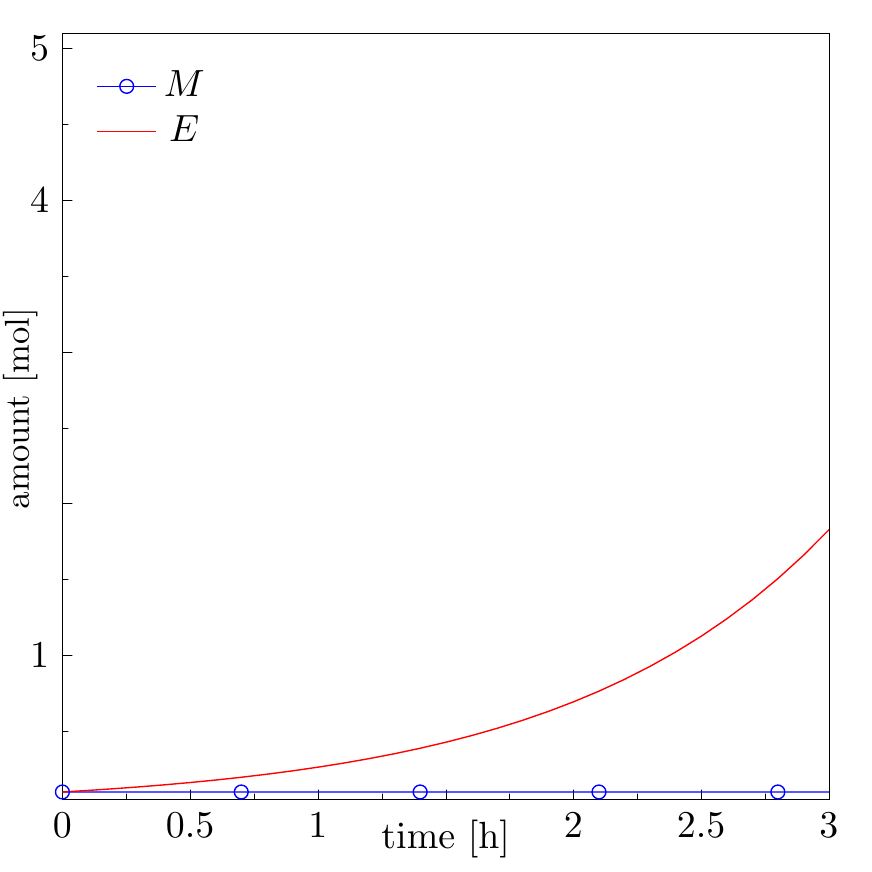}
  \put (50,85) {\textcircled{\smaller B}}
 \end{overpic}
\end{minipage}
\caption{Prediction horizon is set to $p=3.9$ h and step size is $h=0.1$ h.
Plot {\textcircled{\smaller A}} shows the reference solution of the sdeFBA with nominal $k_{\mathrm{cat}}$-values and plot {\textcircled{\smaller B}} shows the robust solution.}\label{fig:robust_growth}
\end{figure}

We visit Example \ref{sec:example_1} a second time to see how the robust solution differs to the nominal one for this minimal example.
We assume all three $k_{\mathrm{cat}}$-values for the reactions include an error up to 20\% of their nominal values.
Hence, we are looking at eight different scenarios consisting of all combinations of extremal values
\begin{align}
 k_{\mathrm{cat}}^j \in \Big\{&  (k_{A,\min},k_{E,\min},k_{M,\min}), (k_{A,\min},k_{E,\min},k_{M,\max}), (k_{A,\min},k_{E,\max},k_{M,\min}), \nonumber\\
			      &  (k_{A,\min},k_{E,\max},k_{M,\max}), (k_{A,\max},k_{E,\min},k_{M,\min}), (k_{A,\max},k_{E,\min},k_{M,\max}), \\
			      &  (k_{A,\max},k_{E,\max},k_{M,\min}), (k_{A,\max},k_{E,\max},k_{M,\max}) \Big\}. \nonumber
\end{align}
We calculated the prediction horizon as $p=3.9~$h via the algorithm described in Section \ref{sec:upper_bound} and set the end-time to $T=3$~h.
As in the previous example, the initial amount of nutrients is chosen so large that they do not limit growth.
Figure \ref{fig:robust_growth} shows three plots for different scenarios. 
The plot {\textcircled{\smaller A}} is the reference solution of the sdeFBA for the nominal values of the catalytic constants.
Due to the increased prediction horizon in comparison to the short end-time $T$, this solution shows only exponential growth.
In particular, we do not observe the artificial linear growth phase (cf. Figure \ref{fig:example_1_results} {\textcircled{\smaller A}}).
Thus, the impact of the fixed end-time $T$ is reduced via the sdeFBA.

The robust solution is presented in plot {\textcircled{\smaller B}}.
Obviously, the inclusion of the uncertainties decreases the overall growth rate significantly.
This is due to the fact that the calculated trajectories are robust.
Thus, they must be feasible for all scenarios at the same time.
The scenario using the minimal values $k_{\mathrm{cat}}^1=\{ k_{A,\min},k_{E,\min}, k_{M,\min} \}$ is in this sense acting as an overall upper bound to the growth rate.
We call this the \emph{minimal scenario}.
The results of the sdeFBA in which we used these minimal values is absolutely identical to the robust solution shown in plot {\textcircled{\smaller B}}.

Thus, we conclude that for such a simple system the robust optimization is unnecessary. 
But we can see impacts of the uncertainties in the Example \ref{sec:core_carbon}, which shows a larger networks including multiple pathways for biomass production.

\subsection{Example: Core Carbon Network}\label{sec:core_carbon}
In this section, we consider a more sophisticated model for the uptake of different nutrients and the transformation of these into biomass.
It corresponds to the core network inside most microorganisms.
This example was previously studied in \cite{waldherr2015} with help of the deFBA under dynamic nutrient conditions.
We only present one environmental setting and compare the results for the different methods.
All transport, metabolic and biomass reactions are shown in the Tables \ref{tb:carbon_reactions}-\ref{tb:biomass_production}.
The different enzymes are either labeled $E_j$ for catalyzing metabolic reactions or $T_j$ for transport reactions.
The catalytic constants $k_{\mathrm{cat}}$ for the enzymes are derived from typical values in metabolism (\cite{bar2011moderately}\cite{schomburg2002brenda}\cite{milo2010bionumbers}).
The constants for the biomass reactions are taken from measurements of translation progression rates inside \textit{E.coli} \cite{young1976}.
The model considers ribosomes R as the essential part of the enzymatic machinery producing all biomass components and itself.
Structural components are collected in the macromolecule S.
Identifying S with the exterior membrane of the cell it also represents the surface area and thus limits the exchange of O2, D, and E.
In this model, the structural component must be more than 35\% of the total biomass for the cell to be working correctly,
\begin{align}
 0.35b^T P \leq S.\label{eq:bcc_example}
\end{align}

For the rdeFBA we assume that three of the $k_{\mathrm{cat}}$-values include a measurement error.
We have chosen the values
\begin{align}
 k_{\mathrm{cat,E_G}} = 1800~1/\mathrm{h}\,d_{\mathrm{E_G}}, ~k_{\mathrm{cat,E_T}} = 1800~1/\mathrm{h}\,d_{\mathrm{E_T}},~k_{\mathrm{cat,R,R}} = 0.2~1/\mathrm{h}\, d_{\mathrm{R,R}},
\end{align}
where $d_{\mathrm{E_G}} \in [1,1.3],d_{\mathrm{E_T}} \in [0.7, 1.3]$ and $d_{\mathrm{R,R}} \in [1,10]$.
Large values for E$_{\mathrm{T}}$ make the production of ATP utilizing O2 more effective, while $\mathrm{E_G}$ is essential for the production of the components G.
The $k_{\mathrm{cat,R,R}}$-value describes the efficiency of self-reproduction of the ribosomes R.
Larger values lead to a quicker growth and deFBA simulations have shown that this may also impact the decision if the cell prefers to stay in exponential growth or goes to a linear growth phase.
We call any phase in which any of the transporter enzymes $\mathrm{T_j}$ are beeing produced an exponential growth phase, while in a linear growth phase only the best yield macromolecule S is produced.
Notice that  the cell will only produce enough S to fulfill the biomass composition constraint \eqref{eq:bcc_example} in an exponential growth phase.
In a linear growth phase production of S is maximized, as the structural component has the largest biomass yield.

This network already captures a large variety of possible behaviors and decisions for the biomass composition.
Varying the environmental conditions, especially the presence of oxygen, leads to very different enzyme levels and pathway usage.
Thus, we are using oxygen dynamics comparable to an aerated batch process
\begin{align}
 \frac{\mathrm{d}}{\mathrm{d}t} \mathrm{O2_{ext}} = V_0 - \gamma_0 \mathrm{O2_{ext}} + V_{\mathrm{O2}},
\end{align}
with the oxygen inflow $V_0 = 20$ mol/min, the ventilation rate $\gamma_0 =0.4$ 1/min, and the cellular uptake flux $V_{\mathrm{O2}}$, which is an optimization variable.

We use the same objective functional as before, namely
\begin{align}
 \underset{Y(\cdot),\Ptil(\cdot),{\tilde{V}}(\cdot)}{\max}~ J( T, Y_0, \Ptil_0) = \underset{Y(\cdot),\Ptil(\cdot),{\tilde{V}}(\cdot)}{\max}~  \int_{0}^{T} b^T \Ptil ~\mathrm{d}t \label{eq:objective_func_e4}.
\end{align}
We set the initial environment to Carb1(0)=2 mol, Carb2(0)=30 mol and $\mathrm{O2_{ext}}(0)=50$ mol.
We use the Resource Balance Analysis \cite{goelzer2011cell} to identify suitable initial values $\Ptil_0 (Y_0)$ leading to an optimal growth rate in this environment for the nominal $k_{\mathrm{cat}}$-values.
This ensures, that differences between the robust solution and the non-robust only originate in the modeling method.
We constrained the amount of macromolecules at time zero to be $b^T P(0)=0.005$ g with the scaling factor for the macromolecules chosen to be $\alpha=100$.
The resulting initial values are shown in Table \ref{tb:biomass_production}.
We calculated the prediction horizon as $p_{\mathrm{up}}\approx 29.9$ min with the algorithm presented in Section \ref{sec:upper_bound} and use it for the sdeFBA as well as the rdeFBA.
The step size for the discretization was chosen as $h=0.5$ min.

\begin{figure}
\hspace*{-2cm}
\begin{minipage}{0.5\paperwidth}
\centering
biomass accumulation\\
\includegraphics[scale=.8]{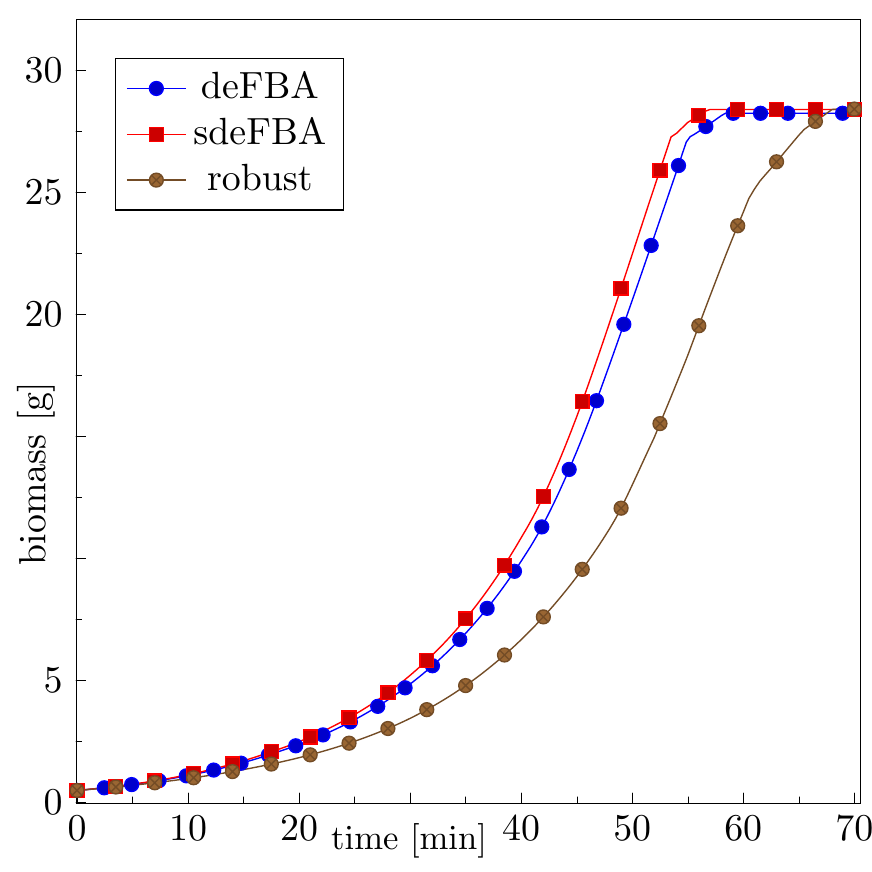}\\
sdeFBA results\\
\includegraphics[scale=0.8]{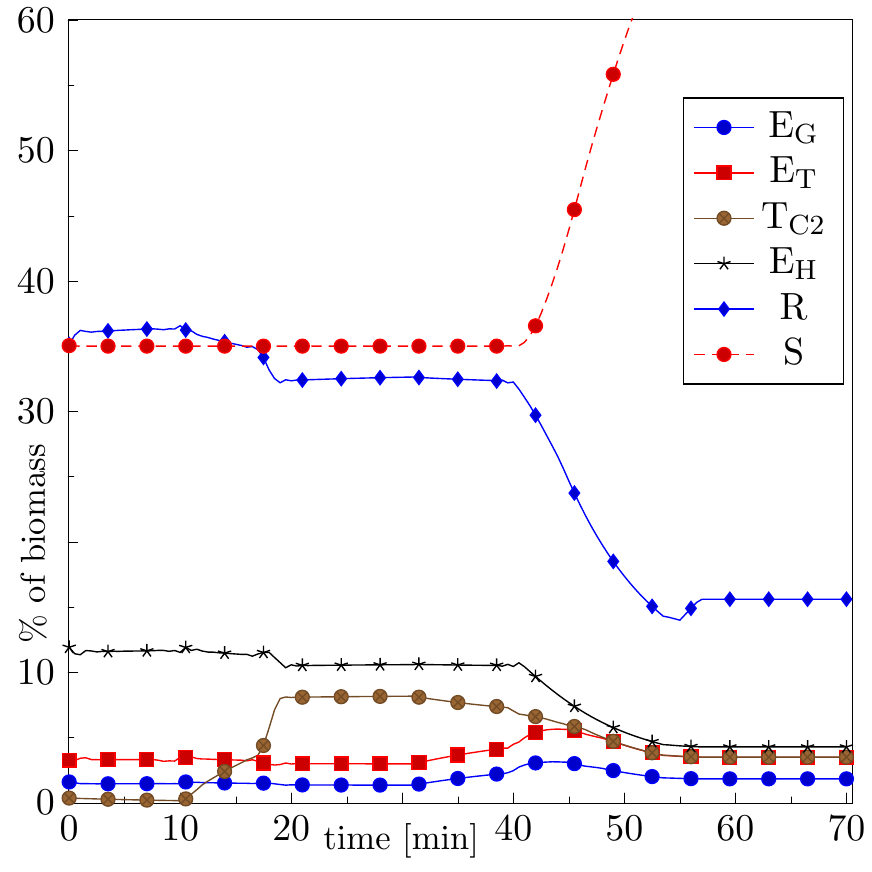}
\end{minipage}\hspace*{-3cm}\begin{minipage}{0.5\paperwidth}
\centering
nutrient dynamics\\
\includegraphics[scale=0.8]{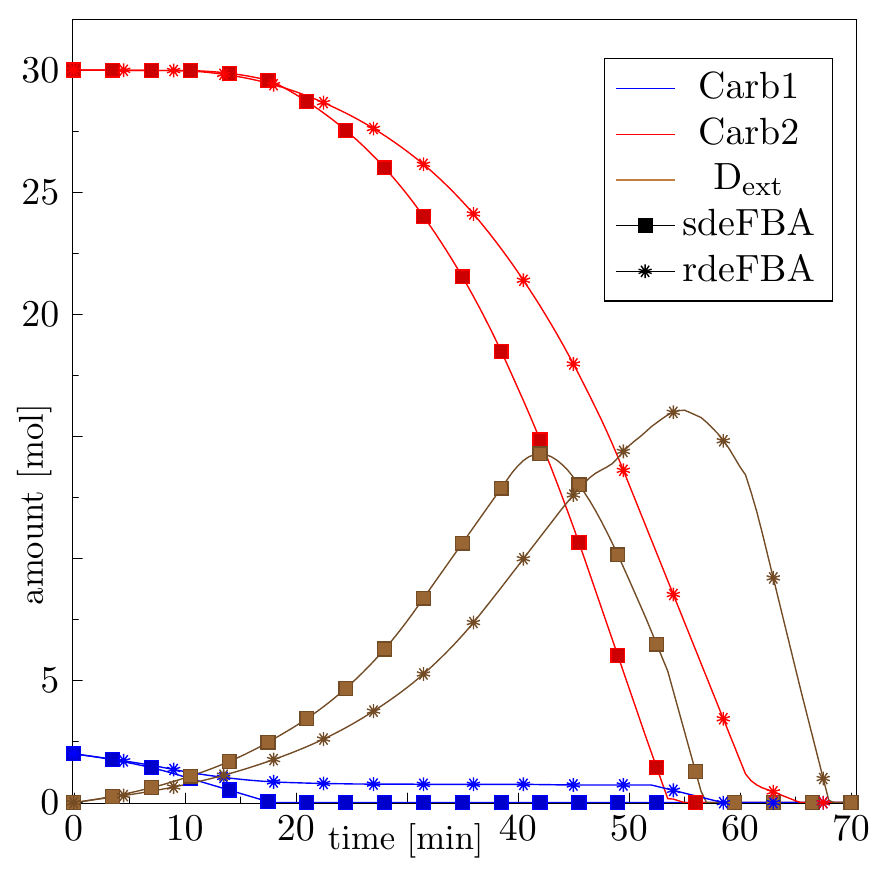}\\
rdeFBA results\\
\includegraphics[scale=0.8]{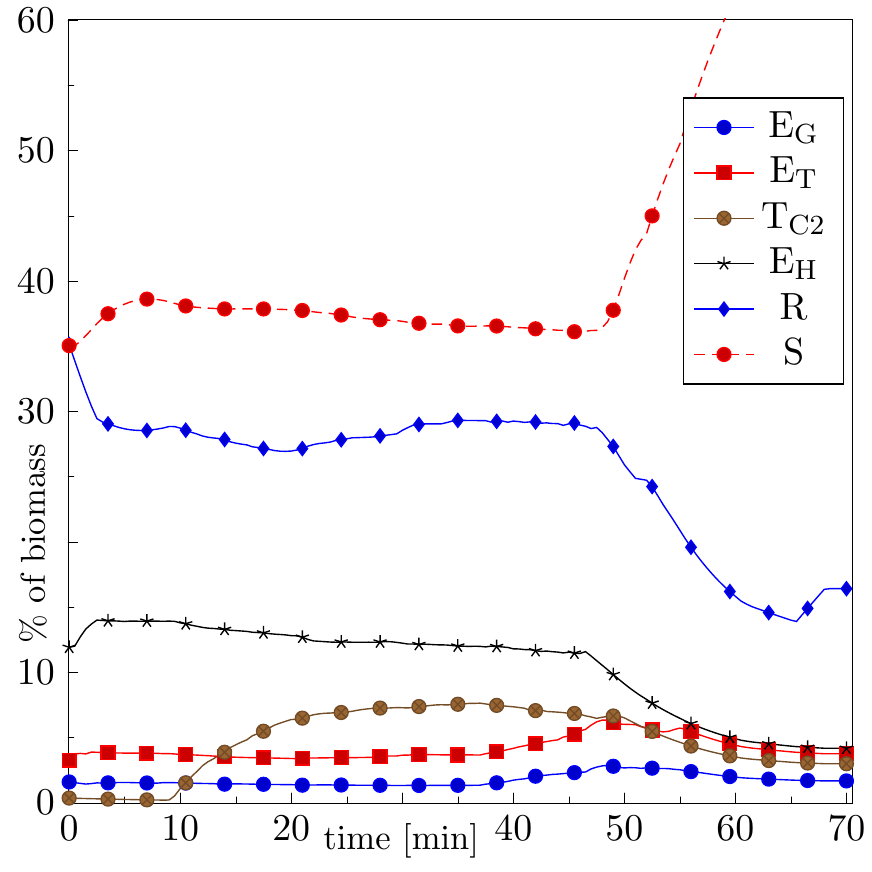}
\end{minipage}
\caption{Comparison of Simulation results for the different methods. (Top left) Plot of biomass $b^T P$ over time for the three methods. (Top right) Comparison of the results for the external nutrients. ($\square$) represents the results for sdeFBA and ($*$) shows the rdeFBA results. (Bottom left) Biomass composition for sdeFBA. (Bottom right) Biomass composition for robust deFBA. }\label{fig:results_example_4}
\end{figure}

The result for the different methods are shown in Figure \ref{fig:results_example_4}.
At the top left the biomass $b^T \Ptil(t)$ is plotted for the different optimizers.
The sdeFBA delivers qualitatively the same results as the deFBA with minor numerical differences.
Hence, we only show the sdeFBA results in more detail, since the enzyme plots for the deFBA are indistinguishable from the ones of sdeFBA.
The overall growth rate of the robust deFBA is much lower than the one observed in the non-robust methods. 
The cause of this becomes clear when looking at the other plots in Figure \ref{fig:results_example_4}.
In the top right plot we compare the time development of the nutrients.
While we can also observe a time delay as in the left plot, we also see that the robust solution does not deplete the Carb1 source at first.
This is very surprising as the parameters are chosen in a fashion to make Carb1 the preferred carbon source.
The two plots at the bottom of Figure \ref{fig:results_example_4} show the time development of the biomass composition and give more insight into the different solutions.
While the sdeFBA predicts an initial increase in Ribosomes along with an increase in enzyme production capacity, the rdeFBA focuses on the production of structure S.
This behavior can be explained with the uncertainty $d_{\mathrm{R,R}}$.
The robust solution optimizes for all scenarios at once and four of those include a 10-fold increased efficiency for the self replication of the ribosomes.
Therefore, the solver decides against an initial production of ribosomes.
Furthermore, the rdeFBA predicts an increased E$_{\mathrm{T}}$ production for increased ATP production from oxygen.  
The sdeFBA solution on the other hand keeps the amount of structural proteins to the enforced minimum \eqref{eq:bcc_example} until the prediction horizon includes the time of nutrient depletion.
Only then the solution predicts an increase in structure molecules S as it is the most effective biomass component in terms of nutrients to biomass conversion.
Overall it is very confirming to see that all solutions reach an identical end-value in total biomass.
Hence, the effectiveness of the production is not affected by our changes in the catalytic constants.
But at the same time we can observe changes in the quality of the solutions and even more in the value of the objective \eqref{eq:objective_func_e4}.

\begin{table}
\caption{Exchange and metabolic reactions, together with their rate constants and catalytic constants $k_{\mathrm{cat}}$. For reversible reactions we use the same value for forward and backward reactions.}\label{tb:carbon_reactions}
\centering
\begin{tabular}{lcc}
\toprule
Reaction & Enzyme & $k_{\mathrm{cat}}~[1/\mathrm{min}]$ \\
\bottomrule
\toprule
Exchange reactions \\
\midrule
Carb1 $\rightarrow$ A & $\mathrm{T}_{\mathrm{C1}}$ & 3000 \\
Carb1 $\rightarrow$ A & $\mathrm{T}_{\mathrm{C2}}$ & 2000 \\
$\mathrm{F}_{\mathrm{ext}} \rightarrow \mathrm{F}$ & $\mathrm{T_F}$ & 3000\\
$\mathrm{H \rightarrow H_{Ext}}$ & $\mathrm{T_H}$ & 3000 \\
$\mathrm{O2_{ext} \rightarrow O2}$ & S & 1000 \\
$\mathrm{D \leftrightarrow D_{Ext}}$ & S & 1000 \\
$\mathrm{E \leftrightarrow E_{Ext}}$ & S & 1000 \\
\bottomrule
\toprule
Metabolic reactions\\
\midrule
A + ATP $\rightarrow$ B & $\mathrm{E_B}$ & 1800\\
B $\rightarrow$ C + 2 ATP + 2 NADH & $\mathrm{E_C}$ & 1800 \\
C $\leftrightarrow$ 2ATP + 3D & $\mathrm{E_D}$ & 1800 \\
C + 4NADH $\leftrightarrow$ 3E & $\mathrm{E_E}$ & 1800 \\
B $\rightarrow$ F & $\mathrm{E_F}$ & 1800 \\
C $\rightarrow$ G & $\mathrm{E_G}$ & 1800 \\
G + ATP + 2NADH $\leftrightarrow$ H & $\mathrm{E_H}$ & 1800 \\
G $\rightarrow$ 0.8C + 2NADH & $\mathrm{E_N}$ & 1800 \\
O2 + NADH $\rightarrow$ ATP & $\mathrm{E_T}$ & 1800 \\
\bottomrule
\end{tabular}
\end{table}

\begin{table}
\caption{List of biomass producing reactions. All reactions are catalyzed by the Ribosomes $R$. Included are the according $k_{\mathrm{cat}}$ values the weights $b$ of the macromolecules.}\label{tb:biomass_production}
\centering
\begin{tabular}{lccc}
\toprule
Biomass reactions & $b~[\mathrm{g / mol}]$ & $k_{\mathrm{cat}}~[1/\mathrm{min}]$ & $P(0)~[\mu\mathrm{mol}] $\\
\bottomrule
\toprule
400 H + 1600 ATP $\rightarrow~\mathrm{T_{C1}}$ & 4 & 2.5 &  23.588\\ 
1500 H + 6000 ATP $\rightarrow~\mathrm{T_{C2}}$ & 15 & 0.67& 0.0\\ 
400 H + 1600 ATP $\rightarrow~\mathrm{T_{F}}$ & 4 & 2.5 &0.0 \\
400 H + 1600 ATP $\rightarrow~\mathrm{T_{H}}$ & 4 & 2.5 & 0.0 \\
500 H + 2000 ATP $\rightarrow~\mathrm{E_B}$ & 5 & 2 & 39.314 \\
500 H + 2000 ATP $\rightarrow~\mathrm{E_C}$ & 5 & 2 & 36.611 \\
1000 H + 4000 ATP $\rightarrow~\mathrm{E_D}$ & 10 & 1 & 19.845  \\
1000 H + 4000 ATP $\rightarrow~\mathrm{E_E}$ & 10 & 1 & 0 \\
2000 H + 8000 ATP $\rightarrow~\mathrm{E_F}$ & 20 & 0.5 & 2.702 \\
500 H + 2000 ATP $\rightarrow~\mathrm{E_G}$ & 5 & 2 &0.0 \\
4000 H + 16000 ATP $\rightarrow~\mathrm{E_H}$ & 40 & 0.25& 14.887  \\
500 H + 2000 ATP $\rightarrow~\mathrm{E_N}$ & 5 & 2 &  0.0 \\
500 H + 2000 ATP $\rightarrow~\mathrm{E_T}$ & 5 & 2 & 13.673 \\
4500 H + 1500 C + 21000 ATP $\rightarrow~\mathrm{R}$ & 60 & 0.2& 30.994 \\
250 H + 250C + 250F + 1500 ATP $\rightarrow~\mathrm{S}$ & 7.5 & 3& 233.333  \\
\bottomrule
\end{tabular}
\end{table}

\section{Discussion}
In this paper, we have introduced two related methods for the simulation of metabolic networks coupled with gene expression.
The proposed short-term deFBA method \eqref{eq:sdeFBA_continous} can be applied to processes in which the end-time of the process is a priori unknown or too small to get a feasible exponential solution.
To decrease the numerical effort in each iteration, we are further investigating the effects of varying prediction horizons.
Basically, we want to find the minimal prediction horizon such that the sdeFBA can produce solution with exponential growth.
We are especially interested in comparing the predictions from the sdeFBA with experiments investigating the actual time span bacteria plan ahead.
The work \cite{zampar2013temporal} observes a decrease in the CO2 production rate of yeast cells 1.5 hours before the main nutrient source, glucose, is depleted.
Yeast must therefore be capable of sensing the nutrient concentration and their gradients.
This shows that the cells are preemptively reacting to environmental changes, which can be identified as the implementation of a prediction horizon inside natural metabolic networks.

The second method we presented is the extension of the sdeFBA to utilizing a robust optimization technique.
We called this the robust deFBA \eqref{eq:rdeFBA}.
It represents a framework to study the impact of uncertainties in metabolic networks and the resulting changes in predicted fluxes.
It inherits the idea of the multi-stage nonlinear model predictive control \cite{lucia2013}.
We have chosen this approach over other robust approaches for its simplicity and easy implementation in combination with the sdeFBA.
Furthermore, the description of uncertainties by discrete events is perfectly suited for our exemplary application with measurement errors in the catalytic constants.
It is very simple to see that is sufficient to use only the maximal deviations from the nominal values of the $k_{\mathrm{cat}}$-values for the creation of the different scenarios.
While already leading to a considerable reduction in the size of the scenario tree, the amount of scenarios still increases two-fold with each additional uncertainty.
Hence, it is only possible to use either a short prediction horizon or a limited number of uncertainties inside the rdeFBA.

In Example \ref{sec:example_3} we have presented an equivalence between the rdeFBA results and the sdeFBA result using the minimal scenario.
This effect can not be reproduced in more complex networks.
We have checked the results of individual sdeFBA problems arising in Example \ref{sec:core_carbon} for all scenarios, whereas, none of these results is identical to the robust solution.
Hence, the effects of the uncertainties can only be observed via the full robust problem including all scenarios at once.

Most importantly, it is obvious that the prediction from the non-robust methods violate the enzyme capacity constraints including the measurement errors.
Therefore, the modeled cells can not reach the growth rate predicted by the non-robust methods.
Additionally, in Example \ref{sec:core_carbon} we have seen that the inclusion of measurement errors led to different usage of the pathways.
This shows that the rdeFBA can serve as a tool to analyze the sensitivity of a model to measurement errors.

\section*{References}
\bibliographystyle{elsarticle-num}
\bibliography{Literatur.bib}
\end{document}